\documentclass[draft,12pt,oneside]{article}
\usepackage{parskip, url, epsf, graphicx, setspace, latexsym, mathrsfs, amsfonts, amssymb, amsmath, wasysym, faktor, geometry, amsthm, xypic, mathptmx}
\usepackage[all]{xy}
\usepackage[mathscr]{euscript}

\title{Genus-Type-Theory}

\newtheorem{thm}{Theorem} 

\newtheorem{prop}[thm]{Proposition}
\newtheorem{cor}[thm]{Corollary}
\newtheorem{definition}[thm]{Definition}
\newtheorem{rmk}{Remark}

\def\On{{\mathbb O}}

\def\Z{{\mathbb Z}}
\def\Q{{\mathbb Q}}

\def\K{{\mathbb K}}

\def\mmp2{{\mathfrak{m}}_p / {\mathfrak{m}}_p^2}
\def\F{{\mathcal{F}}}

\def\f{ \ell }

\def\F{\mathcal{F}}
\def\K{\mathcal{K}}

\def\p{\mathfrak{p}}

\author{John Basias }
\begin{document} 


\maketitle

\section*{Abstract}

We introduce the existence of a Genus-Type Theory that generalizes classical genus theory by linking fractional ideals of number fields to structures built from their Galois groups and associated Diophantine equations, as formally stated in Theorem 30 and Remark 7.

\section*{Introduction }

Let $ K$ be an algebraic number field with Galois closure over $\Q$. We begin by introducing the notion of a $ K$-type Diophantine equation and denote by $\F$ the collection of all such equations. Our goal is to demonstrate the existence of an isomorphism between the group $I_K$ of fractional ideals in $\On_K$ , and a structure derived from $Gal(K / \Q) \oplus \F$ , modulo an appropriate equivalence relation.

Gauss’s work on genera inspires this framework in the theory of binary quadratic forms and the resulting Genus Theory. The central aim of this paper is to establish the existence of a general theory linking Diophantine equations to the ideal class group in the broadest possible context: that of arbitrary algebraic number fields over $\Q$ . However, this work is not constructive—it does not provide explicit methods for realizing such connections.

Accordingly, we refer to the existence of such a relationship as a Genus-Type Theory. In particular, we define a Canonical Genus-Type Theory as a case where one can construct a direct correspondence between a group structure derived from $I_K$ and the set $\F$ of Diophantine equations.

We conclude by illustrating the potential applications of this viewpoint, including a brief reconsideration of classical Genus Theory in the setting of imaginary quadratic fields.

This work is inspired by the ideas of Cox \cite{Cox2013}, and further influenced by results from my thesis \cite{Basias2020}, conducted jointly with Professor Victor Kolyvagin. In that work, we developed several themes originating in Kubota’s study of the structure of biquadratic fields \cite{Kubota T.}, as well as in subsequent generalizations due to Sime \cite{Sime}.

\section{Definitions and Preliminaries}

\begin{definition} Let $K$ be an algebraic number field over $\Q$ that has a Galois group 

and let $G=Gal(K/\Q)$ let $n=[K:\Q]$, $m = 2 \cdot n$, $\On_K$ the ring of integers of $K/\Q$

and $Cl_K$ the class group of $K$. Let $I_K$ and $P_K$ the the Ideals and Principal Ideals in $\On_K$ respectively. \end{definition}

\begin{definition} Let us denote by $\Z [Z_m] = \Z [z_1,z_2,...,z_m]$. \end{definition}

\begin{definition} Let us denote by $o_1, ... , o_n \in \On_K$ a pre-determined basis on $\On_K$

as a $\Z$ module. \end{definition}

\begin{definition} Let us denote by $\On_x = \sum_{i=1}^m a_i \cdot z_i$ where $z_1, ... , z_m$ are variables

and $a_1, ... , a_m \in \On_K$. 

We will denote by $|\On_x|$ the $\On_K$ module over $\Z$ generated by $a_1, ... , a_m$ 

Explicitly $\On_x$ is associated directly to the $m$-tuple $(a_1, a_2, ... , a_m)$ \end{definition}

\begin{definition} Let us denote by $\On_o = \sum_{i=1}^n o_i \cdot z_i$. 

Where for each $i>n$ the variable $z_i$ has coefficient $0$. 

Explicitly $\On_o$ is associated to the $n$-tuple $(o_1, o_2, ... , o_n)$ \end{definition}

\begin{thm} Any Ideal in the ring $\On_K$ is a $\Z$ - Module over $\On_K$ 

of degree at most $m = 2 \cdot n$. \end{thm}

{\bf proof:} It is known that $\On_K$ is a $\Z$ module of degree $n$.

Hence it is the case that any principal ideal in $\On_K$ is a $\Z$ Module over $\On_K$ of degree $n$ 

Now it is also known that any ideal $I$ in $\On_K$ can be generated over $\On_K$ 

by two generators $\alpha_1, \alpha_2$. So we have $I = (\alpha_1) + (\alpha_2)$.

Since each of these principal Ideals are generated by $n$ elements as a 

$\Z$ - Module over $\On_K$ it follows that $I$ must always be a $\Z$ - Module over $\On_K$

of degree at most $m = 2\cdot n$.

\begin{rmk} In what follows next is a treatment of $\On_K$ modules over $\Z$ which have at most $m$ generators. In particular, we make a treatment of Ideals in $\On_K$. Since it is known that an ideal in $\On_K$ is generated by $m$ or less elements over $\Z$. \end{rmk}

\begin{rmk} For the rest of this paper $hom(\Z^m, \Z^m)$ is meant to be the class 

of homomorphisms on the variables $z_1, ... , z_m \mapsto z_1, ... , z_m$. We will express $hom(\Z^m, \Z^m)$ as the monoid group $M_m (\Z)$ .\end{rmk}

\begin{thm} Suppose $\On_x = \sum_{i=1}^m a_i \cdot z_i$ . Likewise let $\On_y = \sum_{i=1}^m b_i \cdot z_i$ 

Explicitly then $\On_x$ is associated to $(a_1, ... , a_m)$ and likewise $\On_y$ is associated to $(b_1, ... , b_m)$. 

Then $|\On_y| \subset |\On_x|$ $\Leftrightarrow$ $\exists$ $h \in M_m (\Z)$ that sends $z_1,...,z_m \mapsto z_1, ... z_m$ 

such that $h(\On_x) = \On_y$.

\end{thm}

{\bf proof:} The forward case $( \Rightarrow )$ is obvious for if $|\On_y| = |h(\On_x)|$ then the $b_i$ 

are clearly each immediately a linear combination of the $a_i$.

So we immediately get that $|\On_y| \subset |\On_x|$. We now treat the reverse $( \Leftarrow)$

Since $|\On_y| \subset |\On_x|$ we have that for each $b_i$, $\exists$ $c_{i1}, c_{i2}, ... c_{im} \in \Z$ 

such that $b_i = \sum_{j=1}^m c_{ij} \cdot a_j$. So now let us write $\On_y = \sum_{i=1}^m b_i \cdot z_i$ as:

$\On_y = \sum_{i=1}^m [\sum_{j=1}^m c_{ij} \cdot a_j \cdot z_i ] = \sum_{j=1}^m a_j \cdot [\sum_{i=1}^m c_{ij} \cdot z_i]$. 

Now we know $\On_x = \sum_{j=1}^m a_j \cdot [z_j]$

So we can make $h \in M_m (\Z) $ as $h(z_j) \mapsto \sum_{i=1}^m c_{ij} \cdot z_i$ 

or explicitly $h = (c_{ij})_{1 \leq i,j \leq m}$ and the result naturally follows.

\begin{cor} For two polynomials $\On_x , \On_y$, we have $|\On_x| = |\On_y|$ $\Leftrightarrow$ 

$\exists$ $\tau_1, \tau_2 \in M_m (\Z)$ such that $\tau_1(\On_x) = \On_y$ and $\tau_2(\On_y) = \On_x$ \end{cor}

\begin{definition} For two linear polynomials $\On_x, \On_y$ we say that $\On_x \sim \On_y$ if

(i) $|\On_x| = |\On_y|$ 

(ii) $\exists$ $\tau_1, \tau_2 \in M_m (\Z)$ such that $\tau_1(\On_x) = \On_y$ and $\tau_2(\On_y) = \On_x$ \end{definition}

\begin{definition} Let us denote by $\K$ all $\On_K$ modules over $\Z$ generated by $m$ or less elements. \end{definition}

\begin{cor} There is an isomorphism between $\K$ and 

$\{\On_x: \On_x=\sum_{i=1}^m a_i \cdot z_i , a_i \in \On_K \}$ under the above equivalence defined in 9. \end{cor}

\begin{definition} Let us denote by $\K^I$ and $\K^P$ all Ideals in $\On_K$ and all principal Ideals in $\On_K$ respectively expressed as $\On_K$ modules over $\Z$. 

We note by Theorem 6 that they are generated by at most $m = 2 \cdot n$ generators. \end{definition}

\begin{definition} We say $f \in \Z[Z_m]$ is of type $K$ if there exists a linear polynomial 

$\On_f = a_1 \cdot z_1 + ... + a_m \cdot z_m$ , $a_1, ... , a_m \in \On_K$ , such that:

$f=N_{K/\Q}(\On_f)$ where $G$ does not act on the variables. 

That is, $\On_f$ is related to the $m$-tuple $(a_1, ... , a_m)$ and $f$ is the Norm-Form of $\On_f$. \end{definition}

\begin{definition} Let us denote by $\F$ the set of all $K$-type polynomials. \end{definition}

\begin{definition}
For any $\On_K$-module over $\Z$, $X$, consider its orbits under the action of $G = Gal(K / \Q)$.

We fix in advance a canonical representative for each orbit and denote it by $X^0$. Then, for each $\sigma \in G$, every element of the orbit can be written as $\sigma(X^0)$, which we denote by $X^\sigma$.

To ensure that these representatives are well defined and correspond one-to-one with the orbits of $X$ under $G$, we impose the following equivalence relation:

If there exist $\sigma_1, \sigma_2 \in G $ such that $X^{\sigma_1}$ and $ X^{\sigma_2}$ lie in the same orbit of $X$, then we declare $ X^{\sigma_1} = X^{\sigma_2}$. In this case, we write $ X^{\sigma_1} \sim X^{\sigma_2}$.

Under this equivalence, the set $\{ X^\sigma : \sigma \in G \} $ corresponds precisely to \\ the set of orbits of $X$.
\end{definition}

\begin{rmk} I have not provided a specific method for predetermining such canonical representatives of $X$, $X^0$ [$\Z$ modules over $\On_K$]. For the purposes of this paper, however, I only need that they are hypothetically pre-determined beforehand. \end{rmk}

\begin{definition}
Let $f$ be a $K$-type Diophantine equation.

We assume that $ f = N_{K / \Q }(\On_f) $ for some linear polynomial $\On_f$.

Let us write $X = |\On_f | $. We then define a specific orbit of $\On_f$, denoted $\On^0_f$, to be the orbit such that $|\On^0_f| = X^0 $, as in Definition 15.

If there are multiple orbits of $\On_f$ with image $X^0$, we simply choose one particular representative.

We define the linear polynomials $\On^{\sigma}_f = \sigma(\On^0_f) $, for each $\sigma \in G $.

Again, as in Definition 15, if it happens that $ |\On^{\sigma_1}_f| = |\On^{\sigma_2}_f| $, then we impose the \\
equivalence: $\On^{\sigma_1}_f \sim \On^{\sigma_2}_f$ .

This restriction is well defined and one-to-one under the equivalences defined in Definitions 14 and 9. Furthermore, under these equivalences, we have:

$\{ \On^{\sigma}_f : \sigma \in G \}$ is equivalent to the set of orbits of $| \On_f | $ .
\end{definition}

\begin{definition} In all of the following sections, we will always be working over the equivalences defined in 14, 9, 15, and 16. \end{definition}

\begin{definition} For a $K$-Type polynomial $f$, that is $f = N_{K/\Q}[\On_f^0]$ we denote the ordered pair $(f,\sigma) = f^{\sigma}$ where $\sigma \in Gal(K / \Q)$ 

We make the association $f^\sigma \leftrightarrow \On_f^{\sigma}$. \end{definition}

\begin{definition} For two elements $f^{\sigma_1}, g^{\sigma_2}$ we make the further equivalence: 

$f^{\sigma_1} \sim g^{\sigma_2}$ if it is the case that: $|{\On_f}^{\sigma_1}| = |{\On_g}^{\sigma_2}| $. \end{definition}

\begin{definition} We denote by $\F_o^I = \{ f^{\sigma} : \sigma \in Gal(K/\Q), f \in \F \}$. \end{definition}

\begin{rmk} We have that $ \F_o^I$ is bijective to $\K$ under the above equivalences. \end{rmk}

\begin{definition} For any $f^{\sigma} \in \F_o^I$ and $\tau \in Gal(K/\Q)$ we define the action $\tau (f^{\sigma}) = f^{\tau \cdot \sigma}$. \end{definition}

\begin{rmk} We have that $ \F_o^I$ is closed under an action of $Gal(K / \Q)$. \end{rmk}

\section{On $K$-type polynomials }

\begin{definition} For a function $h \in M_m (\Z)$ we denote by $h \circ f$ to be:

$h \circ f = N_{K/\Q} [ h(O_f)]$, where $h$ acts on the variables $z_1, ... , z_m$. \end{definition}

\begin{rmk} We can always retrieve $h \circ f$ without factoring $f$ into linnear polynomials. 

For consider $f(z_1, ... , z_m)$ then we have that $h \circ f(z_1,...,z_m) = f(h(z_1),h(z_2),...,h(z_m))$ .\end{rmk}

\begin{definition} Consider a function $h \in M_m (\Z)$ , we say $f$ is a factor of $h \circ f $.

We say this because there is the natural inclusion: $|h(\On_f)| \subset |\On_f|$. \end{definition}

\begin{thm} Consider $|\On_g| \subset |\On_f|$. Let $f = N_{K/\Q}(\On_f), g= N_{K/\Q}(\On_g)$. 

Then $f$ will always be a factor of $g$. \end{thm} 

{\bf proof:} This is clear from the prior section. Let us pick $h$ such that $h(\On_f) = \On_g$

as linear polynomials. The result now follows.

\begin{definition} Let us denote by $\F^I = \{ f^\sigma \in \F^I_o : |\On^{\sigma}_f| \in \K^I \}$ \end{definition}

\section{Multiplication on $\F^I$ and "Genus-Type-Theory" }

\begin{definition} For two elements $f^{\sigma_1}, g^{\sigma_2} \in \F^I$ we define the product of $f^{\sigma_1} * g^{\sigma_2}$ through the following procedure: 

Let us pick two elements $\tau_1, \tau_2 \in M_m (\Z)$ such that $\tau_1(f) = \tau_2( g)= q$. We do this in a certain manner such that this maps to some $q^{\sigma_3} \in \F^I$.

Consider first the function $ q \in \F$. Then $\tau_1$ and $\tau_2$ are chosen again such that $\tau_1(f) = \tau_2( g) = q$ with one additional condition on the multiplication:

We must have that by $|\On^{\sigma_1}_f| \cdot |\On^{\sigma_2}_g|$ lies in the orbit of $|\On^0_q|$.

We note that since $ |\On^{\sigma_1}_f| \cdot |\On^{\sigma_2}_g| \subset |\On^{\sigma_1}_f|$ and $ |\On^{\sigma_1}_f| \cdot |\On^{\sigma_2}_g| \subset |\On^{\sigma_2}_g|$ , 

that it must be the case by Theorem 24 that such $\tau_1, \tau_2$ exist.

So we now define the product to be $q^{\sigma_3}$ which is a representative of the module 

$|\On^{\sigma_1}_f| \cdot |\On^{\sigma_2}_g| = |\On^{\sigma_3}_q|$

As a result, we get that this is represented as $q^{\sigma_3} \in \F^I_o$.

Finally we make the product $f^{\sigma_1} * g^{\sigma_2} = q^{\sigma_3}$

Once More the reasons for this are $f^{\sigma_1} \leftrightarrow \On^{\sigma_1}_f$ and

that $g^{\sigma_2} \leftrightarrow \On^{\sigma_2}_g$ , and lastly that:

$f^{\sigma_1} * g^{\sigma_1}$ we define as $q^{\sigma_3} \leftrightarrow \On^{\sigma_3}_q$, where$|\On^{\sigma_3}_q| = |\On^{\sigma_1}_f| \cdot |\On^{\sigma_2}_g|$

We note that $q^{\sigma_3}$ represents a module that is the product of two ideals.

Hence $|\On^{\sigma_3}_q| \in \K^I$ so indeed: $q^{\sigma_3} \in \F^I \subset \F^I_o$.

\end{definition}

\begin{cor} The structures on $\K^I$ and $\F^I$ are isomorphic up to the equivalences and the multiplication given in Definition 26. \end{cor}

\begin{definition} Let us denote by $\F^P \subset \F^I$ to be the associated set to $\K^P$. \end{definition}

\begin{cor} The multiplication defined on $\F^I$ is closed in $\F^P$ \end{cor}

\begin{thm} [Genus Type Theory] We have $\F^I / \F^P \simeq \K^I /\K^P \simeq Cl_K$ via the prior multiplication defined above. \end{thm}

\begin{definition} Let us denote by $ \F^C = \F^I / \F^P$. \end{definition}

\begin{rmk} $\F^C$ gives the structure of $Cl_K$ via a proper multiplicative equivalence defined on $Gal(K / \Q) \oplus \F$. \end{rmk}

\section{Canonical Genus-Type Theory }

\begin{definition} A Group structure formed from $\F^I$ that is closed under an action of 

$Gal(K/\Q)$, call it $\F_c^C$ is called canonical if it reduces to an equivalence directly on $\F$, 

the set of $K$-type polynomials.

That is explicitly that if for any two $\sigma_1, \sigma_2 \in G$ we have $f^{\sigma_1} \sim f^{\sigma_2} \in \F_c^C$, 

Then we call $\F_c^C$ a Canonical Genus-Type Theory. \end{definition}

We say this because in such cases, it is the case that $Gal(K/\Q)$ plays no role in the 

strucure of $\F_c^C$

As such, if $\F_c^C$ is canonical, its structure relies solely on equivalences on $\F$.

In this section, we provide several settings where quotient groups of $\F^I$ exhibit a Canonical Genus-Type-Theory.

\subsection{On Algebraic Number Fields with a Galois Closure}

In this section, we treat briefly all $K / \Q$ such that $K$ has a Galois closure over $\Q$. 

\begin{definition} Let us denote as $Cl^N_K = N_{K/\Q} [I_K ] / N_{K/\Q} [P_K]$ where $I_K$ are the ideals 

in $\On_K$ and $P_K$ are the principal Ideals in $\On_K$ \end{definition} 

\begin{cor} $Cl^N_K$ can be related to a quotient group of $\F^I$ 

via $Cl^N_K \simeq N_{K/\Q} [\F^I] / N_{K/ \Q}[\F^P]$. \end{cor}

\begin{definition} Let us denote $\F^C_N = N_{K/\Q} [\F^I] / N_{K/ \Q}[\F^P]$ .\end{definition} 

\begin{thm} $\F^C_N$ always is a Canonical Genus-Type-Theory.\end{thm}

{\bf proof:} We choose any arbitrary element of $\F^C_N$ it can be represented as the image of 

an $f^{\sigma_1} \in \F^I$. Now $f^{\sigma_1}$ is also a representative of an ideal $\alpha \in I_K$.

Consider the image of $\alpha$ under the sequence of surjective homomorphisms:

$I_K \overset{\pi_1} \rightarrow N_{K/\Q}[ I_K] \overset{\pi_2} \rightarrow N_{K/\Q} [I_K] / N_{K/ \Q}[P_K] \simeq Cl^N_K$

via $\alpha \overset{\pi_1} \rightarrow N_{K/\Q} [\alpha]$ and $\pi_2$ be the natural mapping. 

It is clear here that for all $\sigma \in Gal(K/\Q)$, $\pi_1[\alpha] = \pi_1[\sigma \alpha ]$

Hence, we have the mapping $\pi_2 \circ \pi_1: I_K \rightarrow Cl^N_K$ has the property that $\alpha \sim \sigma(\alpha)$

for the image of any $\alpha \in I_K$ in $Cl^N_K$.

Now let us consider $f^\sigma \in \F^I$ since we have an equivalence between $\F^I$ and $I_K$ we can

likewise, make the homomorphisms:

$\F^I \overset{\pi'_1} \rightarrow N_{K/\Q}[ \F^I] \overset{\pi'_2} \rightarrow N_{K/\Q} [\F^I] / N_{K/ \Q}[\F^P] \simeq \F^N_C$

Now for the appropriate choice of $f^{\sigma} \in \F^I$, $f^{\sigma} \leftrightarrow \alpha$ , where again $\alpha \in I_K$ 

We have due to our equivalence in the structures that $\forall \tau \in Gal(K/\Q)$ 

$\tau(\alpha) \sim \alpha \in Cl^N_K$ , it is also the case that $\forall \tau \in Gal(K/\Q)$ , $\tau(f^\sigma) \sim f^\sigma \in \F^N_C$.

Now $\tau(f^\sigma) = f^{\tau \circ \sigma}$. Letting $\sigma = \sigma_1$ and $\tau = \sigma_2 \circ {\sigma_1}^{-1}$.

We have shown that for any two $\sigma_1, \sigma_2 \in Gal(K / \Q)$ it is the case that $f^{\sigma_1} \sim f^{\sigma_2} \in \F^C_N$. 

By definition then, $\F^C_N$ is a Canonical Genus-Type-Theory

\subsection{On Kubota-Type Class Groups}

\begin{definition} Let us denote by $L$ a mutli-quadratic extension of $\Q$ that is explicitly 

$\exists$ $d_1, ... d_r \in \Z$ such that $L= \Q[ \sqrt{d_i} : 1 \leq i \leq r ]$. \end{definition}

\begin{definition} Let us denote by $Cl^T_L $ to be the factor group of $Cl_L$ defined as:

$Cl_L / [\prod_{[L:K] = 2} N_{L/K}(Cl_L)]$. We will call this group a Kubota-Type Class Group. 

We say this because the structure of $Cl^T_L$ was inspired by identities introduced

initially by Kubota in \cite{Kubota T.} and then later generalized in my joint work with 

Victor Kolyvagin in Chapter~1, Section~2 of \cite{Basias2020}.

\end{definition}

\begin{definition} Let us call ${\F_t}^C \subset \F^C$ to be the group that represents

the image of $\prod_{[L:K] = 2} N_{L/K}(Cl_L)$ in $\F^C$. \end{definition}

\begin{definition} Let us denote ${\F_T}^C = \F^C/{\F^C}_t$. 

Furthermore, let us denote by ${{\F_T}^C}_{[2]} = {\F_T}^C/ [{\F_T}^C]^2$, i.e. the two-torsion part of ${\F_T}^C$. \end{definition}

\begin{thm} It is the case that ${{\F_T}^C}_{[2]}$ is a Canonical Genus-Type Theory. \end{thm}

{\bf proof:} We show that ${\F_T}^C / {[{\F_T}^C]}^2$ is Canonical.

For any two $\sigma_1, \sigma_2$ consider $\p_1 = f^{\sigma_1}, \p_2 = f^{\sigma_2} \in {{\F_T}^C}_{[2]}$ which are the images of 

the ideals $(\p_1), (\p_2)$ of $\On_L$. We have since the group ${{\F_T}^C}_{[2]}$, is two-primary 

we have that $\p_1 = \p_1 ^{-1}, \p_2 = \p_2^{-1}$ 

Also by the construction of $\F^I$ we have that for any $\sigma, \tau \in G$, $\tau ( f^\sigma ) = f^{\tau \cdot \sigma}$

In the group $\F^I$ and, as a consequence, also for all of its quotient groups. 

So now consider $\sigma_3 \in G$ such that $\sigma_3 \cdot \sigma_1 = \sigma_2$, hence it is the case $\sigma_3 [f^{\sigma_1}] = f^{\sigma_3 \circ \sigma_1} = f^{\sigma_2}$

Hence it must also be the case then that $\sigma_3(\p_1) = \p_2$, since $\p_1 \leftrightarrow f^{\sigma_1}$ and $\p_2 \leftrightarrow f^{\sigma_2}$.

We have that $\p_1 \cdot \p_2 = \p_1 \cdot \sigma_3 (\p_1)$ which is equivalent to $N_{L/K_1}(\p_1)$ 

where $K_1$ is the subfield of $L$ fixed by the orbit of $\sigma_3$.

This lands in the image of ${\F_t}^C$ in the group ${{\F_T}^C}_{[2]}$, which is a part of its kernel. 

As such $\p_1 \cdot \sigma_3 (\p_1) \sim 1$, where we are considering our target group to be ${{\F_T}^C}_{[2]}$.

Now we have that since $\sigma_3 (\p_1) = \p_2$ it follows $\p_1 \cdot \p_2 \sim 1$.

So we have that $\p_2 \sim \p_1^{-1}\sim \p_1 \in {{\F_T}^C}_{[2]}$

Finally we have shown that $\forall \sigma_1, \sigma_2 \in G$ , $f^{\sigma_1} \sim f^{\sigma_2} \in {{\F_T}^C}_{[2]}$.

By the definitions then ${{\F_T}^C}_{[2]}$ is a Canonical Genus-Type Theory.

\begin{rmk} Theorem 41 explicitly shows that a Canonical Genus-Type Theory exists on

any $Cl_k / [Cl_k]^{2}$ when $k$ is a quadratic extension of $\Q$. 

This is so since clearly $N_{k/\Q} (Cl_k) =1$, and as a result we get that:

$Cl^T_k = Cl_k/[ N_{k/\Q} (Cl_k)] = Cl_k$ and as such then $Cl^T_k / [Cl^T_k]^2 = Cl_k/ [Cl_k]^2$ . \end{rmk}

\section{On the Further Development and Application}

Let $f_o = N_{K/\Q} [\On_o]$ as defined in section one. Since every $\On_K$ module over $\Z$ in 

$m$ variables is a subset of $|\On_o|$. We have by Theorems 7 and 24 that:

\begin{thm} We can represent the set $\F = \{ h[f_o] : h \in M_m (\Z) \}$ 

That is explicitly that for any $f \in \F$, $\exists h \in M_m (\Z)$ such that $h[f_o] = f$. \end{thm}

{\bf proof:} Let $\On_f$ be the linear polynomial such that $f = N_{K/\Q}[\On_f]$. 

We have clearly that $\On_f \subset \On_o$.

The result now follows directly from Theorems 7 and 24. 

Now let us assume $\F_c$ is a Canonical-Genus-Type Theory constructed from 

the set of $K$ -type polynomials $\F$. 

\begin{prop} In order to create the Canonical Genus Type Theory, one would need to explicitly complete the
construction of two equivalences. \end{prop} 

{\bf First:} One would have to determine when the two polynomials $f, g \in \F$ 

would be equivalent in $\F_c$.

{\bf Second:} One would have to state the multiplication of two polynomials $f$ and $g$. 

From Theorem 7 and Definition 26, one would have to find the appropriate:

$\tau_1, \tau_2 \in M_m (\Z)$ such that: $\tau_1(f) = \tau_2(g) = f*g = q$

Where $f*g$ represents the product of $f$ and $g$ in the Canonical Genus Type Theory $\F_c$.

\subsection{Classical Genus Theory and Genus-Type Theory}

We will briefly promote some of the concepts and methods of Genus-Theory for imaginary Quadratic Fields through the methods of Genus-Type Theory developed in the previous chapters.

This section is to give the reader an idea of how one can go about explicitly defining the multiplication [in generality] on two $K$-Type Diophantine equations under the proper settings and equivalences, applying the methods stated in Proposition 1.43. 

\begin{definition} Let $D_K$ be the field discriminant of an imaginary quadratic field $K$. \end{definition}

\begin{definition} Let us denote by $d_\f = \f^2 \cdot D_K$ for some pre-determined $\f \in \Z, \f>1$.  \end{definition}

\begin{definition} We will define the Order of Conductor $\f$ in $\On_K$ as $\On_\f = \Z + \f \cdot \On_K$ . \end{definition}

\begin{rmk} For all statements here on forward. We will just refer to $d_\f$ 

as $d$ and assume $\f$ is implied. \end{rmk}

\begin{thm} $\On_\f$ can be expressed as the $\Z$-module in two variables over $\On_K$ 

as $\On_\f = [1, \frac{d-\sqrt{d} }{2}]$. \end{thm}

\begin{definition} We will refer to $\On_\f$ from here on forward as $ \On_\f = x_1 + \frac{d-\sqrt{d}}{2} \cdot x_2$ 

a linear polynomial in two variables over $\On_K$. \end{definition}

\begin{rmk} In the classical sense then the order of conductor $\f$ will be $ |\On_\f|$ .

We make this distinction to make a treatment of the Genus Theory 

through the Genus -Type -Theory. \end{rmk}

The next Theorems 49 and 54 are found in \cite{Cox2013}, Chapter~7, Section~B

\begin{thm} [Structure Theorem] Any proper, invertible ideal of $\On_\f$ can be expressed as $a \cdot \Z + \frac{b+\sqrt{d}}{2} \cdot \Z$ where:

(i) $a \in \Z_{>0}$

(ii) $b^2 \equiv d (\text{mod }4 \cdot a)$ 

(iii) $gcd(a,b, \frac{b^2 - d}{4 \cdot a}) = 1$ .\end{thm}

\begin{prop} Replacing $\Z \rightarrow -\Z$ in the second position we can also say 

all proper invertible ideals are expressible as $a \cdot \Z + \frac{-b- \sqrt{d}}{2} \cdot \Z$. \end{prop}

\begin{definition} For a proper invertible ideal of $\On_\f$ , $\alpha = a \cdot \Z + \frac{-b+\sqrt{d}}{2} \cdot \Z$ we will define 

the linear polynomial $\On_\alpha = a \cdot x_1 + \frac{b-\sqrt{d}}{2} \cdot x_2 $. and in the classical sense then $|\On_\alpha| = \alpha$. \end{definition}

\begin{definition} As in earlier sections, we define all $h \in M_m (\Z)$ to act on the variables of linear polynomials with coefficients in $\On_K$ in m or less variables. 

In particular, we are dealing with lattices, so we restrict $m = 2$. \end{definition}

\begin{definition} Let $h_{\alpha} \in M_2 (\Z)$ to be the element such that $h_{\alpha} (\On_\f) = \On_{\alpha}$. 

We know from section one that such an element exists simply because $\alpha \subset \On_\f$ .

Explicitly $h_{\alpha} = \begin{pmatrix}
a & \frac{b-d}{2} \\
0 & 1 \\
\end{pmatrix} $ . That is explicitly $x_1 \rightarrow a\cdot x_1 + \frac{b-d}{2} \cdot x_2$ , $x_2 \rightarrow x_2$ . \end{definition}

{\bf proof: } $h_{\alpha} ( x_1 + \frac{d-\sqrt{d}}{2} y) = h_{\alpha}(x_1) + (\frac{d-\sqrt{d}}{2}) \cdot h_{\alpha}(x_2) = a \cdot x_1 + \frac{b-d}{2}\cdot x_2 + \frac{d-\sqrt{d}}{2} \cdot x_2$

and one easily sees this is equivalent to $a \cdot x_1 + \frac{b-\sqrt{d}}{2} \cdot x_2 = \On_{\alpha}$ 

\begin{thm} If $a \cdot x^2 + b \cdot xy + c \cdot y^2$ is a primitve quadratic form of discriminant $d$ then $[a, \frac{-b + \sqrt{d}}{2}]$ is a proper, invertible ideal of $\On_\f$. \end{thm}

Furthermore $b^2 - 4 \cdot a \cdot c = d$ so it follows $c=\frac{b^2 - d}{4 \cdot a}$ which leads to the following identity. 

\begin{cor} $N_{K / \Q}(\On_\alpha) = N_{K/\Q} (a \cdot x + \frac{b- \sqrt{d}}{2} \cdot y) = a^2 x^2 + ab \cdot xy + ac \cdot y^2$. Hence:

A proper, invertible ideal of $\On_\f$ represented as $[a, \frac{-b + \sqrt{d}}{2}]$ can be represented as 

the primitive quadratic form expressible as: $\frac{1}{a} \cdot N_{K/\Q}(\On_\alpha) = a \cdot x^2 + b \cdot xy + c \cdot y^2$ .

\end{cor}

We have the next theorem from Cox \cite{Cox2013} Chapter~3, Section~A, which explains in simple terms how to make the product of two primitive quadratic forms of the same discriminant

\begin{thm} for two primitive quadratic forms $a \cdot x^2 + b\cdot xy + c \cdot y^2$ and 

$a' \cdot x^2 + b'\cdot xy + c' \cdot y^2$, we make the product to be the primitive quadratic form:

$aa' \cdot x + B\cdot xy + \frac{B^2 - d}{4aa'} y^2$, where $B$ is the unique number modulo $2aa'$ such that: 

\begin{align}
B \equiv b (mod \text{ } 2a) \\
B \equiv b' (mod \text{ }2a') \\
B^2 \equiv d (mod\text{ } 4aa')
\end{align}

.\end{thm}

Since $aa' \cdot x^2 + B\cdot xy + \frac{B^2 - d}{4aa'} y^2$ is a primitive quadratic form, we have by the corollary that it is represented as the proper invertible ideal:

$\gamma = [aa', \frac{-B+\sqrt{d}}{2}]$ . Let $\alpha = [a, \frac{-b + \sqrt{d}}{2}]$ , $\alpha' = [a', \frac{-b' + \sqrt{d}}{2}]$

A calculation in Cox\cite{Cox2013}, Chapter~7, Section~B shows that $\alpha \cdot \alpha' = \gamma$, this is not hard

to show, and we put the calculation below for the sake of completeness: 

Let $\Delta = \frac{-B+ \sqrt{d}}{2}$ then it is easily seen that $\alpha = [a, \frac{-b + \sqrt{d}}{2}] = [a, \Delta]$ and similarly $\alpha' = [a', \Delta]$ by the equivalences in the theorem above.

As a result $\alpha \cdot \alpha' = [aa', a\Delta , a' \Delta, \Delta^2]$ and we note that $\Delta^2 = - B \cdot \Delta$ hence we get:

$\alpha \cdot \alpha' = [aa', a\Delta , a' \Delta, -B \cdot \Delta]$ and since these forms are primative $gcd(a,a', B) = 1$ 

and hence $\alpha \cdot \alpha' =[aa', \Delta] = \gamma$.

We end this section by incorporating the methods of Proposition 43 to arrive at the identity for the product of two primitive quadratic forms. We note that since we have well-defined $\alpha, \alpha'$ and $\gamma$, the first criterion of the proposition has been satisfied.

We will set $\On_\alpha = a \cdot x + \frac{b - \sqrt{d}}{2} \cdot y$, $\On_{\alpha'} = a' \cdot x + \frac{b' - \sqrt{d}}{2} \cdot y$ and we will retrieve $\On_{\alpha \cdot \alpha'} $ and the primitive quadratic form it expresses in another way. 

We note by theorem 56 that there exist $k_1, k_2 \in \Z$ such that $b + 2 \cdot k_1 \cdot a = b' +2 \cdot k_2 \cdot a' = B$. 

Now, as in Proposition 43, we have successfully stated the equivalence on $ K$-Type polynomials. We proceed now with the second requirement. That is, we can explicitly state our $\tau_1, \tau_2 \in M_2 (\Z) $ such that:

$\tau_1 (\On_\alpha) = \tau_2(\On_{\alpha'}) = \On_{\alpha \cdot \alpha'}$. A quick calculation produces the following theorem:

\begin{thm} $\tau_1 = \begin{pmatrix}
a' & k_1 \\
0 & 1 \\
\end{pmatrix}$. Which is explicitly $x \rightarrow a'\cdot x + k_1 \cdot y$, $y\rightarrow y$. 

Likewise $\tau_2 = \begin{pmatrix}
a & k_2 \\
0 & 1 \\
\end{pmatrix}$. Which is explicitly $x \rightarrow a\cdot x + k_2 \cdot y$, $y \rightarrow y$ 

induce the mappings of $\tau_1( \On_\alpha) = \On_{\alpha \cdot \alpha'}$ and $\tau_2( \On_{\alpha'} ) = \On_{\alpha \cdot \alpha'}$. 

\end{thm}

We now arrive at the identity that yields the primitive quadratic form of the product of $a\cdot x^2 + b\cdot xy + c \cdot y^2$ and $a' \cdot x^2 + b'\cdot xy + c' \cdot y^2$ explicitly through the methods introduced in section 1.3. We note that the primitive polynomial representing the product of the two ideals $\alpha, \alpha'$ by Corollary 55 can be expressed as: 
$\frac{1}{aa'} \cdot N_{K/\Q}( \On_{\alpha \cdot \alpha'} )$

We also note first that by Definition 7:

$N_{K/\Q}(\On_\f) = N_{K/\Q}(x + \frac{d-\sqrt{d}}{2} \cdot y) = x^2 +d \cdot xy + \frac{d^2-d}{4}\cdot y^2$

It is the case $\tau_1 \circ h_\alpha (\On_\f) = \On_{\alpha \cdot \alpha'}$. Explicitly:

$ \frac{1}{aa'} \cdot N_{K/\Q}( \On_{\alpha \cdot \alpha'} ) = \frac{1}{aa'} [ [N_{K/\Q} (\tau_1 \circ h_{\alpha}(\On_f))] = \frac{1}{aa'} [ \tau_1 \circ h_{\alpha} [N_{K/\Q} (\On_\f)]$

We have that $\tau_1 \circ h_\alpha = \begin{pmatrix}
a & \frac{b-d}{2} \\
0 & 1 \\
\end{pmatrix} \cdot \begin{pmatrix}
a' & k_1 \\
0 & 1 \\
\end{pmatrix} = \begin{pmatrix}
a \cdot a' & k_1 \cdot a + \frac{b-d}{2} \\
0 & 1 \\
\end{pmatrix}$

We also know that $B = 2 \cdot a \cdot k_1 + b$ so we may rewrite:

$ \tau_1 \circ h_{\alpha} = \begin{pmatrix}
a \cdot a' & k_1 \cdot a + \frac{b-d}{2} \\
0 & 1 \\
\end{pmatrix} = \begin{pmatrix}
a \cdot a' & \frac{2 \cdot k_1 \cdot a + b-d}{2} \\
0 & 1 \\
\end{pmatrix} = \begin{pmatrix}
a \cdot a' & \frac{B-d}{2} \\
0 & 1 \\
\end{pmatrix}$

We arrive at the explicit identity of the product of the two primitive forms 

setting: $ x \rightarrow {\bf aa'(x) + (\frac{B-d}{2})y }$ and $y \rightarrow y$

\begin{thm} The product of the two primitive quadratic forms 

$a\cdot x^2 + b\cdot xy + c \cdot y^2$ and $a' \cdot x^2 + b'\cdot xy + c' \cdot y^2$ is in fact expressible as:

$\frac{1}{aa'} \cdot[ \tau_1 \circ h_\alpha (x^2 +d\cdot xy + \frac{d^2-d}{4}\cdot y^2 ) ]$, which is explicitly:

$\frac{1}{aa'} \cdot[ {\bf [ aa'(x) + (\frac{B-d}{2})y ]} ^2 + {d \cdot [\bf aa'(x) + (\frac{B-d}{2})y ] } y + \frac{d^2-d}{4}\cdot y^2 ]$.

\end{thm}

\begin{rmk} We should note that in the case of the Classical Genus Theory, the explicit Identity for the product of two Quadratic Forms is known, as exhibited in Theorem 56. 

I've included the method of arriving at the identity through the concepts of Genus-Type-Theory to promote these concepts, as they are potentially applicable in more general settings. \end{rmk}

\clearpage
\addcontentsline{toc}{chapter}{Bibliography}

\nocite{*}
\bibliographystyle{amsplain}

\bigskip
\noindent\textsc{Brooklyn College, CUNY} \\
\noindent\textit{E-mail address}: \texttt{john.basias@brooklyn.cuny.edu}  \\
\noindent\textit{Alternate e-mail}: \texttt{jbasias@gmail.com}

\end{document}